\documentclass[12pt]{article}

\usepackage{a4wide}
\usepackage{amssymb}
\usepackage{amsfonts}
\usepackage{amsmath}
\input xy
\xyoption{arrow}
\xyoption{matrix}

\date{}

\newtheorem{proposition}{Proposition}[section]
\newtheorem{theorem}[proposition]{Theorem}

\newtheorem{corollary}[proposition]{Corollary}

\def\GK{{\rm  GK}\,}
\def\Kdim{{\rm K.dim }\,}

\def\Hom{{\rm Hom}}
\def\der{\partial }

\def\nFM0{{\nu }_{F,M_0}}
\def\nFN0{{\nu }_{F,N_0}}
\def\nGN0{{\nu }_{G,N_0}}

\def\N0{ {\bf N}_0 }

\def\ra{\rightarrow}

\def\Xpm{X^{\pm }}

\def\s{\sigma}
\def\Z{{\mathbb Z }}

\def\l1{{\lambda}_1}

\def\a{\alpha}
\def\a0{ {\alpha }_0}
\def\a1{ {\alpha }_1}

\def\l{\lambda}
\def\o{\omega}

\def\nFGM0{{\nu }_{F,G,M_0}}

\def\nFN0{{\nu}_{F,N_0}}

\def\sm{{\sigma}^m}

\def\sm1{{\sigma}^{-1}}

\def\smtp1{{\sigma}^{-t+1}}

\def\o{\omega }
\def\S1{S^{-1}}

\def\Xpm1{X^{\pm 1}_1}

\def\sPM1{{\sigma }^{\pm 1}}
\def\sMP1{{\sigma }^{\mp 1 }}

\def\b{\beta}
\def\d{\delta}

\def\di{{\rm d.ind}}

\def\L{\Lambda}
\def\O{\Omega}

\def\G{\Gamma}

\def\CA{{\cal A}}

\def\CD{{\cal D}}

\def\Ytm1{Y^{t-1}}
\def\Yim1{Y^{i-1}}

\def\CN{{\cal N}}

\def\CG{{\cal G}}
\def\CH{{\cal H}}

\def\i{{\bf i}}

\def\Der{{\rm Der }}

\def\pad{{\rm pad }}

\def\ker{ {\rm ker } }

\def\D{ \Delta }

\def\SL2Z{ {\rm SL}_2({\bf Z}) }

\def\Gp1{ G^{1 , 1 } }
\def\P11{ P^{-1 , 1 } }
\def\Pp1{ P^{1 , 1 } }

\def\nCLsr{{}^\nu\kern-2pt {\cal L}^{\sigma , \rho  }}
\def\nP{{}^\nu \kern-2pt P}
\def\nL{{}^\nu\kern-2pt L}
\def\nLL{{}^\nu\kern-2pt \Lambda}
\def\nPsr{{}^\nu\kern-2pt P^{\sigma , \rho  }}
\def\nLsr{{}^\nu\kern-2pt L^{\sigma , \rho  }}
\def\nuCL{{}^\nu\kern-2pt  {\cal L}}
\def\nCLsr{{}^\nu\kern-2pt {\cal L}^{\sigma , \rho  }}
\def\nCL1m{{}^\nu\kern-2pt {\cal L}^{-1 , 1  }}
\def\x1nu{x^\frac{1}{\nu}}
\def\xm1nu{x^{-\frac{1}{\nu}}}

\def\CN{{\cal N}}
\def\ra{\rightarrow }

\def\CB{{\cal B}}

\def\CI{{\cal I}}

\def\CC{ {\cal C}}

\def\CH{ {\cal H}}
\def\CP{ {\cal P}}

\def\nAM0{{\nu }_{{\cal A},M_0}}
\def\nAN0{{\nu }_{{\cal A},N_0}}

\def\Kdim{ {\rm Kdim } }

\def\Der{ {\rm Der }}

\def\CP{ {\cal P }}
\def\det{ {\rm det }}

\def\ga{\mathfrak{a}}
\def\gb{\mathfrak{b}}

\def\gn{\mathfrak{n}}

\def\gp{\mathfrak{p}}

\def\j{{\bf j}}

\def\Hom{{\rm Hom}}

\def\JJ{{\bf J}}

\def\di!{\frac{\der^i}{i!}}
\def\dik!{\frac{\der^k_i}{k!}}

\def\gl{\mathfrak{l}}

\def\N{\mathbb{N}}

\def\0{\overline{0}}
\def\1{\overline{1}}

\def\Ln1{\L_{n,\overline{1}}}

\def\oa{\overline{a}}

\def\a1{a_{\overline{1}}}

\def\S{\Sigma}

\def\CU{{\cal U}}

\def\vn1{\overrightarrow{n-1}}

\def\gl{{\rm gl}}
\def\sl{{\rm sl}}

\def\PZ{{\rm PZ}}

\def\mJ{\mathbb{J}}
\def\mI{\mathbb{I}}

\def\ann{{\rm ann}}

\def\ind{{\rm ind}}

\def\K1{{\rm K}_1}

\def\hmI1{\widehat{\mI_1}}
\def\tmI1{\widetilde{\mI_1}}
\def\tmJ1{\widetilde{\mJ_1}}
\def\hB1{\widehat{B_1}}
\def\hCB1{\widehat{\CB_1}}

\def\ga{\mathfrak{a}}

\def\sl2{\mathfrak{sl}_2}
\def\gl2{\mathfrak{gl}_2}

\def\b1{\overline{1}}

\def\PDer{{\rm PDer}}

\def\Sym{{\rm Sym}}

\def\ingp{{\rm in }(\gp )}
\def\outgp{{\rm out }(\gp )}  
\def\vre{\varepsilon}

\begin{document}

\author{V. V. \  Bavula   
}

\title{The Poisson  enveloping algebra and the algebra of Poisson differential operators of a generalized Weyl Poisson algebra}

\maketitle
\begin{abstract}

For a generalized Weyl Poisson algebra $A$, explicit sets of generators and defining relations are presented for  its  Poisson  enveloping algebra $\CU (A)$. Simplicity  criteria are given    for the algebra $\CU (A)$ and algebra  of Poisson  differential operators $P\CD (A)$ on  $A$. The Gelfand-Kirillov dimensions of the algebras $\CU (A)$ and  $P\CD (A)$
 are calculated. It is proven that the algebra $\CU (A)$ is a domain provided that the coefficient ring $D$ of the generalized Weyl Poisson algebra $A$ is a domain of essentially finite type over a perfect field.
 
For the algebra $A$, the set of its minimal primes and the prime radical  are  described  and an equidimensionality criterion is given. 
 For the equidimensional algebra $A$ of essentially finite type, two regularity criteria  are presented. \\

$${\bf Contents}$$
\begin{enumerate}
\item Introduction.
\item Regularity and equidimensionality criteria for commutative generalized Weyl algebras.
\item Generators and defining relations of the Poisson  enveloping algebra of a generalized Weyl Poisson algebra.
\item A simplicity criterion for the Poisson  enveloping algebra of a  generalized Weyl Poisson algebra.
\end{enumerate}

\end{abstract}


\section{Introduction}

In this paper, $K$ is a field of arbitrary characteristic, $P_n=K[x_1, \ldots , x_n]$ is a polynomial algebra in $n$ variables over the field $K$, $D$ is a $K$-algebra and 
$$A=D[X,Y;a]=D[X_1, \ldots , X_N, Y_1, \ldots , Y_N]/(X_1Y_1-a_1, \ldots , X_NY_N-a_N)$$
is a {\bf commutative generalized Weyl algebra} (GWA) of rank $N$ where $X=(X_1, \ldots , X_N)$, $Y=(Y_1, \ldots , Y_N)$ and $a=(a_1,\ldots , a_N)\in D^N$. \\

{\bf The set of minimal primes and an equidimensionality criterion for a GWA $A$.} For an algebra $R$, we denote by $\min (R)$ the set of its minimal primes. Proposition \ref{A13Aug19} gives an explicit description of  the set of minimal primes of a commutative generalized Weyl algebra $A$ via the set of minimal primes of the Noetherian  algebra $D$, it also describes  the prime radical of  the algebra $A$ via the prime radical of the algebra $D$. Proposition \ref{A13Aug19} also gives an equidimensionality criterion for the GWA $A$. \\

{\bf Regularity criteria for the GWA $A$.} In the case when the algebra $D$ is of essentially finite type of pure dimension $d<\infty$ over a perfect field $K$, two regularity criteria are given for the GWA $A= D[X,Y; a]$ of rank $N$, Theorem \ref{14Aug19} and Theorem \ref{13Aug19}. The first criterion, Theorem \ref{14Aug19}, is given via an explicit ideal of the algebra $D$. The second  criterion, Theorem \ref{13Aug19}, is given in terms of a new concept introduced in the paper - a {\em homological sequence} of algebra (every regular sequence is a homological sequence but not the other way round, in general). Proposition \ref{8Aug19} collects some of the properties of homological sequences. \\

{\bf The generalized Weyl Poisson algebras and their Poisson simplicity criterion.}  In \cite{Pois-GWA-Com}, a new large class of Poisson algebras - the class of  generalized Weyl Poisson algebras -  is introduced.

 An associative  commutative algebra $A$ is called a {\bf Poisson algebra} if it is a Lie algebra $(A, \{ \cdot, \cdot \})$ that satisfies the {\em Leibniz's rule}: For all elements $a,x,y\in D$, 
 $$\{ a, xy\}= \{ a, x\}y+x\{ a, y\}.$$
 The set  $\PZ (A):=\{ a\in A\, | \, \{a, x\}=0$ for all $x\in A\}$ is called the {\bf Poisson centre} of $A$, it is a subalgebra of $A$. Let $\Der_K(A)$ be the set of $K$-derivations of the associative algebra $A$. Then   
 $$\PDer_K(D):=\{ \d \in \Der_K(D)\, | \, \d (\{ a,b\})= \{ \d (a),b\}+\{ a,\d (b)\} \;\; {\rm  for\; all}\;\; a,b\in D\}$$
 is  the {\em set of derivations}  of the Poisson algebra $D$. Elements of $\PDer (A)$ are called {\bf Poisson derivations} of $A$.\\

 {\it Definition, \cite{Pois-GWA-Com}.} Let $D$ be a Poisson algebra, $\der = (\der_1, \ldots , \der_n)\in \PDer_K(D)^n$ be an $n$-tuple of commuting derivations of the Poisson algebra $D$, $a=(a_1, \ldots , a_n)\in \PZ (D)^n$ 
  be such that $\der_i (a_j)=0$ for all $i\neq j$. The commutative generalized Weyl algebra $$A= D[ X, Y;  a]=
 D[X_1, \ldots , X_n , Y_1, \ldots , Y_n]/(X_1Y_1-a_1, \ldots , X_nY_n-a_n)$$ admits a Poisson structure which is an extension of the Poisson structure on $D$ and is given by the rule: For all $i,j=1, \ldots , n$ and $d\in D$,

\begin{equation}\label{PGWAR1}
\{ Y_i, d\}=\der_i(d)Y_i, \;\; \{ X_i, d\}=-\der_i(d)X_i \;\; {\rm and}\;\; \{ Y_i, X_i\} = \der_i (a_i),
\end{equation}
\begin{equation}\label{PGWAR2}
\{ X_i, X_j\}=\{ X_i, Y_j\}=\{ Y_i, Y_j\} =0\;\; {\rm for\; all}\;\; i\neq j.
\end{equation}
The Poisson algebra is denoted by $A =D[ X, Y; a, \der \}$ and is called the {\bf generalized Weyl Poisson algebra} of degree/rank $n$ (or GWPA, for short)  where $X=(X_1, \ldots , X_n)$ and $Y= (Y_1, \ldots , Y_n)$.\\

Examples of GWPAs are given in Section \ref{PEAGWPA}. A Poisson simplicity criterion is given in \cite[Theorem 1.1]{Pois-GWA-Com}, see Theorem \ref{10Apr16}. This result is used in several theorems of this paper.\\

{\bf Generators and defining relations of the Poisson  enveloping algebra $\CU (\CA )$ of a GWPA $A$.} For an arbitrary Lie algebra $\CG$, every $\CG$-module is a module over its {\bf universal enveloping algebra} $U(\CG )$, which is an associative algebra,  and vice versa. Furthermore, the defining relations of the algebra $U(\CG )$ are precisely the defining relations of the Lie algebra $\CG$ where the Lie bracket $[\cdot , \cdot ]_\CG$ of $\CG$ is `replaced' by the algebra commutator: If $\{ X_i\}_{i\in I}$ is a $K$-basis of the Lie algebra $\CG$   then the algebra $U(\CG )$ is generated by the elements $\{ X_i\}_{i\in I}$ that subject to the defining relations: $$X_iX_j-X_jX_i=[X_i, X_j]_\CG\;\; {\rm for \; all}\;\; i,j\in I.$$ 
Similarly, for each Poisson algebra $\CP$ there is a concept of Poisson module over $\CP$, and every Poisson module over $\CP$ is a module over, the so-called, Poisson  enveloping algebra $\CU (\CP )$, which is an associative algebra, and vice versa. In \cite{GenDefRel-PUEA}, for each Poisson algebra $\CP$ that is defined by generators abd defining relations, explicit sets of generators and defining relations are given for its Poisson  enveloping algebra $\CU (\CP )$, \cite[Theorem 2.2]{GenDefRel-PUEA}, see Theorem \ref{23Jun19}. Using this result, for each GWPA $A$, Theorem \ref{E23Jun19} gives explicit sets of generators and defining relations for the algebra $\CU (A)$.  Examples are considered in Proposition \ref{F23Jun19}, Corollary \ref{aF23Jun19} and Corollary \ref{bF23Jun19}. \\

{\bf The  Poincar\'{e}-Birkhof-Witt Theorem for Poisson algebras.} The (classical) Poincar\'{e}-Birkhof-Witt Theorem states that for each Lie algebra $\CG$ there is a natural isomorphism of   graded algebras,  $${\rm gr}\, U(\CG ) \simeq \Sym (\CG ),$$ where ${\rm gr}\, U(\CG )$ is   the associated graded algebra   of the   universal enveloping algebra  $ U(\CG )$ of the Lie algebra $\CG$ and $\Sym (\CG )$  is the symmetric algebra of $\CG$. For a {\em smooth} Poisson algebra $\CP$,  a
similar result holds \cite[Theorem 3.1]{Rinehart-1963}: 
 $${\rm gr}\, \CU(\CP ) \simeq \Sym_\CP (\O_\CP ),$$ where ${\rm gr}\, \CU(\CP )$ is   the associated graded algebra   of the   Poisson  enveloping algebra  $ \CU(\CP )$ of the Poisson algebra $\CP$ and $\Sym_\CP (\O_\CP )$  is the symmetric algebra of the $\CP$-module  $\O_\CP$ of K\"{a}hler differentials of the associative algebra $\CP$. In fact,  \cite[Theorem 3.1]{Rinehart-1963} holds in sightly more general situation, namely, for the universal enveloping algebra of the Lie-Reinhatr algebra. The pair $(\CP, \O_\CP )$ is an example of a Lie-Reinhart algebra and its universal enveloping algebra $V(\CP, \O_\CP )$ is isomorphic to the Poisson enveloping  algebra (PEA) $\CU (\CP )$, \cite[Section 2, p.197]{Rinehart-1963} (see also \cite{Huebschmann-90}).
 Recently,  it was proven  that the PBW Theorem holds for certain {\em singular} Poisson hypersurfaces,  \cite[Theorem 3.7]{Lambre-Ospel-Vanhaecke-2019}.
 \cite[Theorem 3.5]{GenDefRel-PUEA} 
  states that the  Poincar\'{e}-Birkhof-Witt Theorem holds for all  Poisson algebras (over an arbitrary field). \\

In \cite{Oh-1999}, the Poisson enveloping algebra of a Poisson algebra was introduced as a  universal object in a certain category and an alternative  to Reinhart's proof of its existence was given. For certain classes of Poisson algebras explicit descriptions  of their Poisson enveloping algebras were presented in \cite{Oh-Park-Shin-2002, Umirbaev-2012, Yang-Yao-Ye-2013,Lu-Wang-Zhuang-2015, Lu-Oh-Wang-Yu-2018, Lambre-Ospel-Vanhaecke-2019}. \\

{\bf Simplicity criterion for the algebra $P\CD (A)$ of Poisson differential operators on  the generalized Weyl Poisson algebra $A$.}  Corollary \ref{GWPA10Apr16}
 is a simplicity criterion  for the algebra  of Poisson differential operators  $P\CD (A )$ on  a generalized Weyl Poisson algebra $A$.\\

 {\bf Simplicity criterion for the Poisson enveloping  algebra $\CU (A )$ of the generalized Weyl Poisson algebra $A$.}  Corollary \ref{GWPAY22Mar19} 
 is a simplicity criterion  for the algebra  $\CU (A )$  where $A$ is   a generalized Weyl Poisson algebra.
 If in addition, the coefficient ring $D$ of the generalized Weyl Poisson algebra 
 $A=D[X,Y;a, \der\}$ is a regular domain of essentially finite type over a field of characteristic zero, Theorem \ref{GWPA7Aug19} is an explicit simplicity criterion for the algebra $\CU (A)$. Theorem \ref{GWPA7Aug19}.(4) is a very efficient tool in proving/disproving simplicity of the algebra $\CU (A)$. \\

{\bf The Gelfand-Kirillov dimension of the algebras $\CU (A )$, $P\CD (A)$,    ${\rm gr}\, \CU (A )$ and $ \Sym_A (\O_A )$.} 
 Let $A=D[X,Y;a, \der\}$  be a generalized Weyl Poisson algebra where $D$ is a  domain of essentially finite type over a perfect field.  Let $\O_A$ be the $A$-module of  K\"{a}hler differential of the algebra $A$ and  $ \Sym_A (\O_A )$ be its symmetric algebra.  Corollary \ref{GAAA29Jul19} 
 (resp., Corollary \ref{GBB26Jul19}) gives an exact figure for the Gelfand-Kirillov dimension of the algebras  $\CU (A )$,    ${\rm gr}\, \CU (A )$ and $ \Sym_A (\O_A )$ (resp., $P\CD (A)$, char$(K)=0$).\\
 
  {\bf The algebra $\CU (A )$ is a domain when $A$ is a regular domain of essentially finite type.}  Let $A=D[X,Y;a, \der\}$  be a generalized Weyl Poisson algebra where $D$ is a  domain of essentially finite type over a perfect field. Corollary  \ref{G29Jul19} states that the algebra $\CU (A )$ is a domain provided the algebra $A$ is a  regular domain.


\section{Regularity and equidimensionality criteria for commutative generalized Weyl algebras}\label{REGCRIT}

The aim of this section is to give regularity criteria for a commutative generalized Weyl algebra $A=D[X,Y;a]$ of rank $N$ where $D$ is a commutative algebra of essentially finite type of pure dimension $d<\infty$ over a perfect field $K$. The first criterion, Theorem \ref{14Aug19}, is given in terms of an explicit ideal of the ring $D$. The second  criterion, Theorem \ref{13Aug19}, is more conceptual and is given via a new concept - a homological sequence of algebra. One of the key results in proving Theorem \ref{13Aug19} is Proposition \ref{A13Aug19} which describes the set of minimal primes of the algebra $A$ via the set of minimal primes of the algebra $D$ and gives an equidimensionality criterion for the algebra $A$. \\

{\bf Generalized Weyl algebras, \cite{Bav-GWA-FA-91, Bav-SimGWA-1992, Bav-GWArep}.} Let $D$ be a ring, $\sigma=(\sigma_1,...,\sigma_n)$ be  an $n$-tuple of
commuting automorphisms of $D$,  $a=(a_1,...,a_n)$ be  an $n$-tuple  of elements of  the centre
$Z(D)$  of $D$ such that $\sigma_i(a_j)=a_j$ for all $i\neq j$. The {\bf generalized Weyl algebra} $A=D[X, Y; \sigma,a]$ (briefly GWA) of degree/rank   $n$  is  a  ring  generated
by $D$  and    $2n$ indeterminates $X_1,...,X_n, Y_1,...,Y_n$
subject to the defining relations:
$$Y_iX_i=a_i,\;\; X_iY_i=\sigma_i(a_i),\;\; X_id=\sigma_i(d)X_i,\;\; Y_id=\sigma_i^{-1}(d)Y_i\;\; (d \in D),$$
$$[X_i,X_j]=[X_i,Y_j]=[Y_i,Y_j]=0, \;\; {\rm for \; all}\;\; i\neq j,$$
where $[x, y]=xy-yx$. We say that  $a$  and $\sigma $ are the  sets  of
{\it defining } elements and automorphisms of the GWA $A$, respectively.

 The $n$'th {\bf Weyl algebra},
 $A_n=A_n(K)$ over a field (a ring) $K$ is  an associative
 $K$-algebra generated by  $2n$ elements
 $X_1, ..., X_n,Y_1,..., Y_n$, subject to the relations:
$$[Y_i,X_i]=\d_{ij}\;\; {\rm and}\;\;  [X_i,X_j]=[Y_i,Y_j]=0\;\; {\rm for\,all}\; i,j, $$
where $\d_{ij}$ is the Kronecker delta function.
The Weyl algebra $A_n$ is a generalized Weyl algebra
 $A=D[X, Y; \s ;a]$ of rank $n$ where
$D=K[H_1,...,H_n]$ is a polynomial ring   in $n$ variables with
 coefficients in $K$, $\s = (\s_1, \ldots , \s_n)$ where $\s_i(H_j)=H_j-\delta_{ij}$ and
 $a=(H_1, \ldots , H_n)$.  The map
$$A_n\ra A,\;\; X_i\mapsto  X_i,\;\; Y_i \mapsto  Y_i,\;\;  i=1,\ldots ,n,$$
is an algebra  isomorphism (notice that $Y_iX_i\mapsto H_i$).

It is an experimental fact that many quantum algebras of small Gelfand-Kirillov dimension are GWAs (eg, $U({\rm sl}_2)$, $U_q({\rm sl}_2)$, the quantum Weyl algebra, the quantum plane, the Heisenberg algebra and its quantum analogues, the quantum sphere,  and many others).

The GWA-construction turns out to be a useful one. Using it for large classes of algebras (including the mention ones above) all the simple modules were classified, explicit formulae were found for the  global
and Krull dimensions, their elements were classified in the sense of Dixmier, \cite{Dix}, etc.\\

{\bf A $\Z^n$-grading on a GWPA.} The GWPA of rank $n$,
\begin{equation}\label{AZngr}
A:= D[X,Y; a, \der \} =\bigoplus_{\alpha\in \Z^n}A_\alpha ,
\end{equation}
is a $\Z^n$-graded Poisson algebra where $A_\alpha = Dv_\alpha$,  $v_\alpha =\prod_{i=1}^n v_{\alpha_i}(i)$ and $$v_j(i)=\begin{cases}
X_i^j& \text{if }j>0,\\
1& \text{if }j=0,\\
Y_i^{|j|}& \text{if }j<0.\\
\end{cases}$$

So, $A_\alpha A_\beta \subseteq A_{\alpha +\beta}$ and  $\{ A_\alpha , A_\beta \}\subseteq A_{\alpha +\beta}$ for all elements $\alpha , \beta \in \Z^n$.\\

{\bf The set of minimal primes of a generalized Weyl algebra.} For an algebra $R$, $\min (R)$ is the set of its minimal primes and  $\gn_R$ is its prime radical. We say that a commutative algebra $R$ has {\bf pure dimension} $d$ if the Krull dimension of all factor algebras $R/ \gp$, where $\gp \in \min (R)$, is $d$. We also say the the algebra $R$ is {\bf equidimensional}. Let $\CC_R$ be the set of all regular elements of the algebra $R$, (i.e. the set of non-zero-divisors).

 Poposition \ref{A13Aug19} gives an explicit description of  the set of minimal primes of a commutative generalized Weyl algebra $A$ via the set of minimal primes of the Noethirian  algebra $D$, it also describes  the prime radical of  the algebra $A$ via the prime radical of the algebra $D$.

\begin{proposition}\label{A13Aug19}
 Let $D$ be a commutative Noetherian ring, $A=D[X,Y; a]$ be a commutative GWA of rank $N$, for each $\gp \in \min (D)$, let 
 $$\ingp =\{ j\, | \, a_j\in \gp, 1\leq j \leq N\}\;\; {\rm and }\;\;
 \outgp =\{ j\, | \, a_j\not\in \gp, 1\leq j \leq N\}.$$ Then 
 \begin{enumerate}
\item $\min (A)=\coprod_{\gp \in \min (D)} \min (A, \gp )$ where 
\begin{eqnarray*}
\min (A, \gp )&:=&\{ P\in \min (A)\, | \, \gp \subseteq P\} = \{ P_{\gp, \varepsilon} \, | \, \vre =(\vre_1, \ldots ,\vre_{|\ingp |})\in \{ \pm \}^{{|\ingp |}}  \},\\
P_{\gp, \varepsilon}&:=& \begin{cases}
(\gp , v_{\vre_1}(i_1), \ldots , v_{\vre_{|\ingp |}}(i_{|\ingp |}))& \text{if } \ingp \neq\emptyset,\\
(\gp )& \text{if } \ingp =\emptyset ,\\
\end{cases}
\end{eqnarray*}
$v_+(j)=X_j$ and $v_-(j)=Y_j$, 
\begin{enumerate}
\item for all $P\in \min (A, \gp )$, $P\cap D = \gp$,
\item $|\min (A, \gp )|=2^{|\ingp |}$ and $\min (A)|\geq |\min (D)|$,  
\item $\min (A)|= |\min (D)|$ iff $\ingp = \emptyset$ for all $\gp \in \min (D)$ iff $\min (A) = \{ A\gp \, | \, \gp \in \min (D)\}$.
\end{enumerate}
\item $\gn_A= A\gn_D$ where $\gn_A$ and $\gn_D$ are the prime radicals of the rings $A$ and $D$, respectively. In particular, the ring $A$ is reduced iff the ring $D$ is so.
\item For each minimal prime $\gp \in \min (D)$, the set $\min (A, \gp )$ is the set of all the minimal primes over the ideal $A\gp$ of $A$, and the ideal $A\gp = \bigcap_{P\in \min (A, \gp )}P$ is a semiprime ideal of $A$.
\item $A/P_{\gp, \vre}=D/\gp [v_{-\vre_1}(i_1), \ldots , v_{-\vre_{|\ingp |}}(i_{|\ingp |})][(X_{j_1}, \ldots, X_{j_{|\outgp |}}), (Y_{j_1}, \ldots, Y_{j_{|\outgp |}}); $ $ a=(a_{j_1}, \ldots, a_{j_{|\outgp |}}) ]$ is a commutative GWA of rank $|\outgp |$ with coefficients in the polynomial algebra $$D/\gp [v_{-\vre_1}(i_1), \ldots , v_{-\vre_{|\ingp |}}(i_{|\ingp |})]$$ in $|\ingp |$ variables with coefficients in $D/\gp $. In particular, the field of fractions of the algebra $A/P_{\gp, \vre}$ is  $$Q(A/P_{\gp, \vre})=Q(D/\gp ) (X_1, \ldots , X_N)$$ which is the field of rational functions in the variables $X_1, \ldots, X_N$ over the field of fractions $Q(D/\gp )$ of the algebra $D/\gp $, and ${\rm tr.deg}_K Q(A/P_{\gp, \vre })={\rm tr.deg}_K Q(D/\gp) +N$.
\item If, in addition, the algebra $D$  has pure dimension $d$ then the algebra $A$ has pure dimension $d+N$, and vice versa.
\end{enumerate}
\end{proposition}

{\it Proof}. 4. Statement 4 is obvious (use the $\Z^N$-grading of the GWA $A$).

1. (i) {\em All ideals  $P_{\gp, \vre }$ of $A$ are distinc prime ideals}: This follows from statement 4.

(ii) $P_{\gp, \vre }\cap D=\gp$: This follows from the fact that the ideal $P_{\gp, \vre }$ is a homogeneous ideal of the GWA $A$ (with respect of the $\Z^N$-grading) and definitions of the sets $\ingp$ and $\outgp$.

(iii) {\em The set $\min (A)$ contains precisely all the ideals $P_{\gp, \vre }$}: Let $P$ be a prime ideal of $A$. Then the intersection $D\cap P$ is a prime ideal of the algebra $D$. Hence, $\gp \subseteq D\cap P$ for some $\gp \in \min (D)$. Let $\ingp =\{ i_1, \ldots , i_l\}$ and $\outgp =\{ j_1, \ldots , j_m\}$. Then the factor algebra
$$A/A\gp = D/\gp [X_\outgp , Y_\outgp ; a_\outgp ][X_\ingp , Y_\ingp ; 0]$$
 where $X_\outgp =(X_{j_1}, \ldots , X_{j_m})$, $Y_\outgp =(Y_{j_1}, \ldots , Y_{j_m})$,  $X_\ingp =(X_{i_1}, \ldots , X_{i_l})$, $Y_\ingp =(Y_{i_1}, \ldots , Y_{i_l})$ and $a_\outgp = (\overline{a_{j_1}}, \ldots , \overline{a_{j_m}})$ where $\overline{a_\nu }= a_\nu +\gp$. Since 
 $$X_{i_1}Y_{i_1}=0, \ldots , X_{i_l}Y_{i_l}=0,$$ we must have $P_{\gp , \vre}\subseteq P$ for some $\vre$. So, every prime ideal of $A$ contains some prime ideal $P_{\gp , \vre}$  and all the ideals $P_{\gp , \vre}$ are distinct. Now, the statement (iii) follows, and as a result we have statement 1. 
 
 3. By statement 1, $A\gp = \bigcap_{P\in \min (A, \gp )}P$, and statement 3 follows.

2. Statement 2 follows from statements 1 and 3:
\begin{eqnarray*}
A\gn_D &=&A\Big(\bigcap_{\gp \in\min  (D)}\gp \Big)=\bigcap_{\gp \in \min (D)} A\gp \stackrel{{\rm st}.3}{=}\bigcap_{\gp \in \min (D)}\bigcap_{P \in \min (A, \gp )}P\\
 &\stackrel{{\rm st}.1}{=}& \bigcap_{P \in \min (A)}P=\gn_A.
\end{eqnarray*}
The rest is obvious.

5. Statement 5 follows from statement 4: For each minimal prime $P_{\gp , \vre}$ of $A$, 
$$ A/P_{\gp , \vre}=D/\gp [X_\outgp , Y_\outgp ; a_\outgp ][v_{-\vre_1}(i_1), \ldots , v_{-\vre_l}(i_l)]$$
is a polynomial algebra in the variables $v_{-\vre_1}(i_1), \ldots , v_{-\vre_l}(i_l)$ withe coefficients in the GWA $B=D/\gp [X_\outgp , Y_\outgp ; a_\outgp ]$ of degree $m=|\outgp |$ and all the coordinates of the vector $a_\outgp =( a_{j_1}, \ldots , a_{j_m})$ are nonzero (hence regular) elements of the domain $D/ \gp$. Hence, $Q(A/P_{\gp , \vre})=Q(D/\gp ) (X_1, \ldots , X_m)$. Now, 
$$\Kdim (A/P_{\gp , \vre})=\Kdim (B)+l.$$
{\sc Claim.} $\Kdim (B) = \Kdim (D/\gp ) + m = d+m$.\\
 
By the Claim, $$\Kdim (A/P_{\gp , \vre})=d+m+l=d+n,$$ as required (since $m+l=n$). 

To prove the Claim we use induction on $m$. The result is obvious if $m=0$. Suppose that $m>0$ and that the Claim is true for all $m'<m$. Then 
$$B=B_{m-1}[X_m, Y_m; a_m]=B_{m-1}[X_m, Y_m]/(X_mY_m-a_m).$$
The polynomial algebra $B_{m-1}[X_m, Y_m]$ is a Noetherian domain of Krull dimension  $\Kdim (B_{m-1})+2$. By the Krull's Principal Ideal Theorem, the algebra $B$ is a Noetherian algebra of pure dimension $$\Kdim (B_{m-1})+2-1=\Kdim (B_{m-1})+1= d+m-1+1=d+m$$ 
 since, by induction, $\Kdim (B_{m-1})= d+m-1$, and the Claim follows. $\Box $\\

{\bf Homological sequences.} 

{\it Definition.} Let $D$ be a regular algebra of essentiallu finite type. A sequence of elements $a_1, \ldots , a_N$ in $D$ is called a {\bf homological sequence of length} $N$ if the following two conditions hold:

1. the element $a_1$ is a {\bf homological sequence of length} 1, that is (by definition) the element $a_1$ is a regular element of $D$ such that the factor algebra $D/(a_1)$ is regular provided the element $a_1$ is not a unit, and

2. the images $\oa_2, \ldots , \oa_N$ of the elements $a_2, \ldots , a_N$ in the ring $$D_1:=\begin{cases}
D& \text{if } a_1\in D^\times ,\\
D/(a_1)& \text{if } a_1\not\in D^\times \cup \{ 0\} ,\\
\end{cases}$$
is a homological sequence of length $N-1$. The algebra $D_1$ is called {\bf the algebra of the element $a_1$} in the homological sequence $s$: $a_1, \ldots , a_N$ in $D$. Similarly, the algebra $D_2$ of the element $\oa_2$ in the homological sequence $\oa_2, \ldots , \oa_N$ in $D_1$ is called {\bf the algebra of the element $a_2$ in the homological sequence} $s$ in $D$. Finally, {\bf the algebra $D_N$ of the element $a_N$} in the homological sequence $s$ is defined. The sequence of algebras $D_1, \ldots , D_N$ is  called {\bf the sequence of algebras of the homological sequence} $s$. \\

Notice that the algebras $D_1, \ldots , D_N$ are regular algebras of essentially finite type and $\oa_i\in \CC_{D_{i-1}}$ for all $i=1, \ldots , N$.  Let 
$$\CU (s) :=\{ a_i \, | \, \oa_i \in D_{i-1}^\times \}\;\; {\rm  and}\;\;\CN (s) := \{ a_i \, | \, \oa_i \not\in D_{i-1}^\times \}.$$ The set $\CU (s)$  (resp., $\CN (s)$) is called the {\em set of relative units} (resp., {\em regular non-units}) of $s$. If $\CN (s) = \{ a_{i_1}, \ldots , a_{i_\nu}\}$ where $1\leq i_1<\cdots <a_{i_\nu}\leq N$ then 
\begin{equation}\label{Dinu}
D_{i_\mu}=D/(a_{i_1}, \ldots , a_{i_\nu})\;\; {\rm for \; all}\;\; \mu = 1, \ldots , \nu .
\end{equation}
Furthermore, there is a chain of equalities and epimorphisms  :
$$ D_0=D_1=\cdots = D_{i_1-1}\ra D_{i_1}= \cdots = D_{i_2-1}\ra \cdots D_{i\mu}=\cdots = D_{i_\mu -1}\ra \cdots \ra D_{i_\nu} = \cdots = D_N$$ where the arrows are epimorphisms (eg, $D_{i_\mu -1}\ra D_{i_\mu}\simeq D_{i_\mu -1} / (a_\mu ))$ and $D_0:=D$. The elements in $\CN (s)$ are a regular sequence of the algebra $D$. The algebra  $D_i$ has pure dimension $d-u_i$ where $u_i=\# \{ j\, | \, a_j\in \CN (s), j\leq i\}$.\\

{\bf Regularity criteria for a GWA $A$.} Theorem \ref{14Aug19} is a first regularity criterion for the a GWA $A$. 

\begin{theorem}\label{14Aug19}
Let $D=S^{-1}(P_n/I)$ be a commutative algebra of essentially finite type of pure dimension $d$ $(=n-r$ where $P_n=K[x_1, \ldots , x_n]$, $I=(f_1, \ldots , f_m)$ and $r=r(\frac{\der f_i}{\der x_j})$) over a perfect field $K$ and $A=D[X,Y; a]$ be a GWA of rank $N$. Then the algebra $A$ is regular iff $\sum_{I}\gb_I=D$ where the sum is taken over all subsets $I\subseteq \{ 1, \ldots , N\}$ with $|I|\geq N-d$, $CI=\{ 1, \ldots , N\} \backslash I$ and $\gb_I:= \ga_{CI} a_I$ is an ideal of $D$ where $a_I:= \prod_{i\in I} a_i$,  $a_\emptyset :=1$, and the ideal
 $\ga_{CI}$ of $D$  is generated by all the minors of size $r+|CI |$ of the $(m+|CI |)\times n$ matrix with coefficients in $D$ where the $m+|CI |$ rows are ${\rm grad} (f_1), \ldots ,  {\rm grad} (f_m), {\rm grad} (a_{j_1}), \ldots ,  {\rm grad} (a_{j_{|CI |}})$ where $CI = \{ j_1<\cdots <j_{|CI |}\}$. 
\end{theorem}

{\it Proof}.  By Theorem \ref{A13Aug19}.(5), the algebra $A$ is an algebra of essentially finite type of pure dimension $d+N$ over a perfect field $K$. By the Jacobian Criterion of Regularity, the algebra $A$ is regular iff the Jacobian ideal $\ga = \ga (A)_{d+N}$ of $A$ is equal to $A$. The strategy of the proof of the theorem is as follows. First, we find explicit set of generators of the ideal $\ga$. Then we show that $D\cap \ga = \sum_I\gb_I$, and the theorem  follows. Notice that 
$$ A=S^{-1}K[x_1, \ldots , x_n, X_1, Y_1, \ldots , X_N, Y_N]/(f_1, \ldots , f_m, a_1-X_1Y_1, \ldots , a_N-X_NY_N).$$
The Jacobian matrix of $A$ is of size $(m+N)\times (n+2N)$ with $m+N$ rows as follows 
\begin{eqnarray*}
({\rm grad} (f_1),0,\ldots , 0),  &\ldots ,&   ({\rm grad} (f_m),0,\ldots , 0),  \\
({\rm grad} (a_1), -X_1, -Y_1, 0, \ldots , 0), &\ldots ,&   ({\rm grad} (a_N),  0, \ldots , 0, -X_N, -Y_N).
\end{eqnarray*}  
It follows that 
$$ \ga = \sum_{I, \vre} \ga_{CI} v_{I, \vre}$$
where $v_{I, \vre} = \prod_{i\in I} v_{\vre_i}(i)$, $\vre = (\vre_i)_{i\in I}\in \{ \pm \}^I$ and $v_+(i)=X_i$ and $v_-1(i) = Y_i$ and the set $I$ runs through all the subsets of $\{ 1, \ldots , N\}$ such that $r+N-|I|\leq n$ $(\Leftrightarrow |I|\geq N-d$). The ideal $\ga$ is a sum of homogeneous ideals w.r.t. the $\Z^N$-grading of the GWA $A$. So, 
$$ 1\in \ga \;\; {\rm iff}\;\; 1\in \ga\cap D=\sum_{I, \vre } \ga_{CI} v_{I, \vre} v_{I, -\vre}= \sum_I\ga_{CI} a_I,$$
as required. $\Box$

Let $S_N$ be the symmetric group, i.e. the group of all permutations of the set $\{ 1, \ldots , N\}$. Theorem \ref{13Aug19} is the second regularity criterion for GWAs. 

\begin{theorem}\label{13Aug19}
Let  $D$ be a commutative  algebra of essentially finite type of pure dimension $d$ over a perfect field $K$ and $A=D[X,Y; a]$ be a GWA of rank $N$. Then the following statements are equivalent:
\begin{enumerate}
\item The algebra $A$ is a regular algebra. 
\item  The algebra $D$ is a regular algebra and  the sequence $a_1, \ldots , a_N$ is a homological sequence in $D$ of length $N$. 
\item The algebra $D$ is a regular algebra and every permutation of the sequence $a_1, \ldots , a_N$ is a homological sequence in $D$ of length $N$. 
\end{enumerate}
If one of the equivalent conditions hold then the algebras $D_1, \ldots , D_N$ of the homological sequence $a_1, \ldots , a_N$ are commutative regular equidimensional algebras of essentially finite type. Furthermore, the pure dimension of the algebra $D_i$ is $d-u_i$  where $u_i=\# \{ j \,  | \, a_j\in \CN (s) , 1\leq j \leq i\}$. 
\end{theorem}

{\it Proof}. $(1\Leftrightarrow 2)$ To prove that the equivalence holds we use induction on $N$. Suppose that $N=1$. By Theorem \ref{14Aug19}, the algebra $A$ is regular iff $\sum_I\gb_I=D$ where the sum is taken over all the subsets $I$ of the set $\{ 1\}$ such that $|I|\geq 1-d$. There are two cases to consider: $d=0$ and $d>0$. 

Suppose that $d=0$. Then $|I|\geq 1$, i.e. $I=\{ 1\}$. Now, $D= \sum_{I=\{ 1\} }  \gb_I = \ga_ra_1$ iff $\ga_r=D$ and $a_1\in D^\times$. 

Suppose that $d>0$. Then $|I|\geq 1-1=0$, and so $I=\{ 1\}$ or $I=\emptyset$. 
 Now, $$D= \sum_{I\in \{ \emptyset , \{ 1\}\} } \gb_I = \ga_\emptyset a_1+\ga_{\{ 1\} }a_\emptyset = \ga_ra_1+\ga_{\{ 1\} } \subseteq \ga_r\;\; ({\rm since}\;\; a_\emptyset = 1)$$ iff $\ga_r=D$ and $D= (a_1, \ga_{ \{ 1\} })$ iff the algebra $D$ is regular and the algebra $D/ (a_1)$is regular since $(a_1, \ga_{\{ 1\} })$ is the Jacobian ideal of the $D/(a_1)$ of essentially finite type of pure dimension $d-1$ (by Krull's Ptincipal Ideal Theorem). The proof of the equivalence when $N=1$ is complete. 
 
 Suppose that $N>1$, and the equivalence holds for all nartural numbers $N'<N$. Notice that $A= A_{N-1}[X_1, Y_1; a_1]$ is a GWA of rank 1 where $A_{N-1} = D[X', Y'; a']$ is a GWA of degree $N-1$ where $X'=(X_2, \ldots , X_N)$,  $Y'=(Y_2, \ldots , Y_N)$ and $a' = (a_2, \ldots , a_N)$. Recall that the algebra $D$ is an algebra of essentially finite type of pure dimension $d$. By Proposition \ref{A13Aug19}.(5), the algebra $A_{N-1}$ is an algebra of essentially finite type of pure dimension $d+N-1$. By the case $N=1$, $A$ is regular iff the algebra $A_{N-1}$ is regular, $a_1\in \CC_{A_{N-1}}$ $(\Leftrightarrow a_1\in \CC_D)$ and if $a_1\not\in A_{N-1}^\times$  $(\Leftrightarrow a_1\not\in D^\rtimes)$ then the algebra $$A_{N-1}/(a_1)\simeq D/(a_1)[X', Y' ; \oa']$$ is regular where $\oa'=(a_2+(a_1), \ldots , a_N+(a_1))$ iff $a_1$ is a homological sequence in $D$ and $\oa_2, \ldots , \oa_N$ is a homological sequence in $D_1$. This means that the sequence $a_1, \ldots , a_N$ is  homological sequence in $D$. 
 
 $(2\Leftrightarrow 3)$ For any permutation $\s \in S_N$, $A\simeq D[X^\s, Y^\s ; a^\s ]$ where $X^\s = (X_{\s (1)}, \ldots , X_{\s (N)})$ and $a^\s = (a_{\s (1)}, \ldots , a_{\s (N)})$, and the equivalence follows from the fact that 
 $(1\Leftrightarrow 2)$.  $\Box$

We say that ideals $\ga$ and $\gb$ of an algebra $R$ are {\em co-prime} if $\ga + \gb =R$.

\begin{corollary}\label{a16Aug19}
Let $D$ be a ring of essentially finite type of pure dimension $d$ over a perfect field $K$ and $A=D[X,Y; a]$ be a commutative GWA of degree $N$. The following statements are equivalent.
\begin{enumerate}
\item The algebra $A$ is a regular algebra.
\item $A\simeq \prod_{P\in \min (A)}A/P$ and the algebras $A/P$ are regular.
\item The algebra $D$ is regular, elements in  $\min (A)$ 
are pairwise co-prime, and for each $P=P_{\gp, \vre}$, the GWA 
$$D/\gp [(X_{j_1}, \ldots, X_{j_{|\outgp |}}), (Y_{j_1}, \ldots, Y_{j_{|\outgp |}});  a=(a_{j_1}, \ldots, a_{j_{|\outgp |}}) ]$$ 
is a regular algebra.
\end{enumerate}
\end{corollary}

{\it Proof}. 1. $(1\Leftrightarrow 2)$ The equivalence is true for all algebras of essentially finite type. 

$(2\Leftrightarrow 3)$ The algebra $D$ is reduced ($\gn_D=0$) iff the algebra $A$ is so, by Proposition \ref{A13Aug19}.(2), i.e. $0=\gn_A=\bigcap_{P\in \min (A)}P$. 
 So, $A\simeq \prod_{P\in \min (A)}A/P$ iff the minimal primes of the algebra $A$ are pairwise co-prime. Now, 
 the equivalence $(2\Leftrightarrow 3)$ follows from Proposition \ref{A13Aug19}.(4). $\Box $

\begin{proposition}\label{8Aug19}
Let  $D$ be a commutative regular algebra of essentially finite type of pure dimension $d$ over a perfect field $K$ and $a_1, \ldots , a_N$ be a homological sequence in $D$ of length $N$. Then 
\begin{enumerate}
\item Every permutation of the sequence $a_1, \ldots , a_N$ is a homological sequence in $D$ of length $N$. 
\item  There is a permutation, say $a_{i_1}, \ldots , a_{i_N}$,  of the 
sequence $a_1, \ldots , a_N$ such that 
\begin{enumerate}
\item the algebras $D_0:=D, D_1, \ldots , D_\nu$ are regular where $D_i=D/ (a_1, \ldots , a_i)$ for $i=1,\ldots , \nu$, 
\item $\oa_i\in \CC_{D_i}\backslash D_i^\times$ for $i=1, \ldots , \nu$, and
\item $\oa_{\nu +1}, \ldots , \oa_N\in D_s^\times$ where $\oa_j = a_j+(a_1, \ldots , a_\nu)$.
\end{enumerate}
\end{enumerate}
\end{proposition}

{\it Proof}. 1. Statement 1 follows from Theorem \ref{13Aug19}.

2. Let $D_1, \ldots , D_\nu$ be as in (\ref{Dinu}). Then $D_\nu = D/(a_{i_1}, \ldots , a_{i_\nu})$ where $1\leq i_1<\cdots <i_\nu\leq N$ and $a_j+(a_{i_1}, \ldots , a_{i_\nu})\in D_\nu^\times$ for all $j\in \{ j_1, \ldots , j_{n-\nu }\} = \{ 1, \ldots , N\} \backslash \{ i_1, \ldots , i_\nu \}$, by Theorem \ref{13Aug19}. Then the sequence $a_{i_1}, \ldots , a_{i_\nu }, a_{j_1}, \ldots , a_{j_{n-\nu}}$ is required permutation of the original sequence. $\Box$


\section{ Generators and defining relations of the Poisson  enveloping algebra of a generalized Weyl Poisson algebra}\label{PEAGWPA}

The aim of this section is for an arbitrary generalized Weyl Poisson algebra to give explicit sets of generators and defining relations for its Poisson  enveloping algebra (Theorem \ref{E23Jun19}). We also consider some interesting examples (Proposition \ref{F23Jun19} and Corollary \ref{aF23Jun19}). At the begining of the section, examples  of generalized Weyl Poisson algebras are considered. \\

{\bf Examples of GWPAs.} 1. If $D$ is a  algebra with {\bf trivial} Poisson bracket $( \{ \cdot , \cdot \}=0)$ then any choice of elements $a=(a_1, \ldots , a_n)$ and $\der = (\der_1, \ldots , \der_n)\in \Der_K(D)^n$ such that $\der_i(a_j) = 0$  for all $i\neq j$ determines a GWPA $D[X,Y;a, \der\}$ of rank $n$ where $\Der_K(D)$ is the set of $K$-derivations of the algebra $D$. If, in addition, $n=1$ then there is no restriction on $a_1$ and $\der_1$.\\

2. The {\bf classical Poisson polynomial algebra} $P_{2n} = K[X_1, \ldots , X_n , Y_1, \ldots , Y_n]$ ($\{ Y_i, X_j\} =\d_{ij}$ (the Kronecker delta) and $\{X_i, X_j\} = \{ X_i, Y_j\} = \{ Y_i, Y_j\} =0$ for all $i\neq j$) is a GWPA
\begin{equation}\label{P2nKH}
P_{2n}=K[H_1, \ldots , H_n][X,Y; a, \der \}
\end{equation}
where $K[H_1, \ldots , H_n]$ is a Poisson polynomial algebra with trivial Poisson bracket, $a=(H_1, \ldots , H_n)$, $\der = (\der_1, \ldots , \der_n)$ and $\der_i= \frac{\der}{\der_{H_i}}$ (via the  isomorphism of Poisson algebras $P_{2n}\ra K[H_1, \ldots , H_n][X,Y; a, \der \}$, $X_i\mapsto X_i$, $Y_i\mapsto Y_i$).\\

3. $A= D[X,Y; a,\der \}$ where $D=K[H_1, \ldots , H_n]$ is a Poisson polynomial algebra with trivial Poisson bracket, $a= (a_1, \ldots , a_n)\in K[H_1]\times \cdots \times K[H_n]$, $\der = (\der_1, \ldots , \der_n)$ where  $\der_i=b_i\der_{H_i}$ (where $\der_{H_i}=\frac{\der}{\der H_i}$)  and  $b_i\in K[H_i]$. In particular, $D[X,Y; (H_1, \ldots , H_n), (\der_{H_1}, \ldots , \der_{H_n})\}=P_{2n}$ is the classical Poisson polynomial algebra.

Let $S$ be a multiplicative set of $D$. Then $S^{-1}A\simeq (S^{-1} D)[X,Y; a, \der \}$ is a GWPA. In particular, for $S=\{ H^\alpha \, | \, \alpha \in \Z^n\}$ we have $K[H_1^{\pm 1}, \ldots , H_n^{\pm 1}][X,Y; a, \der \}$. In the case $n=1$, the Poisson algebra $$K[H_1^{\pm 1}][X_1, Y_1; a_1, -H_1\frac{d}{dH_1}\}$$ where $a_1\in K[H_1^{\pm 1}]$  is,  in fact, {\em isomorphic} to  a Poisson algebra   in the paper of Cho and Oh \cite{CHO-OH} which is obtained as a quantization of a certain GWA with respect to the quantum parameter $q$.  In \cite[Theorem 3.7]{CHO-OH}, a Poisson simplicity criterion  is given for this Poisson algebra. \\

{\bf Poisson $\CP$-modules.} Let  $(\CP , \{ \cdot, \cdot \} )$ be a  Poisson algebra. A  left $\CP$-module $M$  is called a  {\bf Poisson   $\CP$-module} if there is a bilinear map
$$ \CP \times M\ra M, \;\; (a, m) \mapsto \d_am$$
which is called a {\em Poisson action} of $\CP$ on $M$ such that for all elements $a,b\in \CP$ and $m\in M$,\\

(PM1) $\d_{\{ a,b\} }= [\d_a, \d_b ]$, \\

(PM2)  $[\d_a, b]= \{ a,b\}$, and \\

(PM3) $\d_{a b}= a\d_b+b\d_a$.\\

{\bf Example.} Every Poisson algebra $\CP$ is a Poisson $\CP$-module where $\d_a= \pad_a$. \\

{\bf The Poisson  enveloping algebra $\CU (\CP )$ of a Poisson  algebra $\CP$.}\\

For each Poisson algebra $\CP$ there is a unique (up to isomorphism) {\em associative} algebra $\CU (\CP )$ such that \\

$\bullet$ {\em every Poisson $\CP$-module is a $\CU (\CP)$-module, and vice versa.}\\

The algebra $\CU (\CP )$ is called the {\bf Poisson  enveloping algebra} (PEA)  of the Poisson  algebra $\CP$.\\

{\bf  Generators and defining relations of the Poisson  enveloping algebra of a generalized Weyl Poisson algebra.}
 For a Poisson algebra $\CP$ which is defined by generators and defining relations (as an associative algebra), Theorem \ref{23Jun19} gives explicit sets of generators and defining relations for the Poisson    enveloping algebra $\CU (\CP )$. 

\begin{theorem}\label{23Jun19}
(\cite[Theorem 2.2.(2)]{GenDefRel-PUEA})   Let $\CP$ be a Poisson algebra, $U (\CP )$ be its  enveloping algebra (as a Lie algebra)  and  $\CU (\CP )$ be its Poisson    enveloping algebra. Then 
\begin{enumerate}
\item $\CU (\CP ) \simeq \CP \rtimes_\pad U (\CP )/\CI (\CP )$ where $\CI (\CP )= (\d_{ab} -a\d_b - b\d_a )_{a,b\in \CP }$ is the ideal of the algebra $\CP \rtimes_\pad U (\CP )$ generated by the set $ \{ \d_{ab} -a\d_b - b\d_a \, | \, a,b\in \CP \}$.
\item If $\CP = S^{-1}K[x_i]_{i\in \L} / (f_s)_{s\in \G}$ where $S$ is a multiplicative subset  of the polynomial algebra $K[x_i]_{i\in \L}$ ($\L$ and $\G$ are index sets). Then the algebra $\CU (\CP )$ is generated by the algebra $\CP$ and the elements $\{ \d_i:= \d_{x_i}\, | \, i\in \L \}$ subject to the defining relations (a)--(c): For all elements $ i,j\in \L$ such that $i\neq j$ and $s\in \G$,
\begin{enumerate}
\item $[\d_i , \d_j] = \sum_{k\in \L} \frac{\der \{ x_i, x_j\}}{\der x_k}\d_k$, 
\item $[\d_i, x_j]= \{ x_i, x_j\}$, and 
\item $\sum_{i\in \L} \frac{\der f_s}{\der x_i}\d_i=0$. 
\end{enumerate}
So, the algebra $\CU (\CP )$ is generated by the algebra $\CP$ and the set $\d_\CP =\{ \d_a\, | \, a\in \CP \}$ subject to the defining relations: For all elenents $a,b\in \CP$ and $\l , \mu \in K$, 
\begin{enumerate}
\item $[\d_a, \d_b ] = \d_{\{ a,b\} }$, 
\item $[\d_a , b] = \{ a,b\}$, 
\item $\d_{ab} = a\d_b + b\d_a$,  
\item $\d_{\l a+\mu b} = \l \d_a+\mu \d_b$ and $\d_1=0$. 
\end{enumerate}
\item The map $\pi_\CP : \CU (\CP )\ra \CD (\CP )$, $a\mapsto a$, $ \d_b \mapsto \pad_b= \{ b, \cdot \}$ is an algebra homomorphism where $a, b\in \CP$ and its image is the algebra $P\CD (\CP )$ of Poisson differential operators of the Poisson algebra $\CP$. 
\item The algebra $\CP$ is a subalgebra of $\CU (\CP )$. Furthermore, $\CU (\CP ) = \CP \oplus \ann_{\CU (\CP )}(1)$ is a direct sum of left $\CP$-modules where $\ann_{\CU (\CP )}(1)=\sum_{i\in \L}\CU (\CP) \d_i$ is the annihilator of the identity element of the Poisson $\CP$-module $\CP$. The Poisson $\CP$-module structure on the Poisson algebra $\CP$ is obtained from the $\CD (\CP)$-module structure on $\CP$ by restriction of scalars via $\pi_\CP$. 
\end{enumerate}
\end{theorem}

\begin{theorem}\label{E23Jun19}
 Let $A= D[X,Y; a, \der \}$ be a GWPA of degree $n$ and $D=S^{-1}K[t_\nu ]_{\nu \in I}/(f_s)_{s\in S}$. Then  the Poisson  enveloping algebra $\CU (A)$ of the Poisson algebra $A$ is generated by the algebra $A$ and elements $\d_{X_1}, \ldots, \d_{X_n}, \d_{Y_1}, \ldots, \d_{Y_n}, \d_{t_\nu}$ (where $\nu \in I$) subject to the defining relations given in statements 1-3: For  all elements $i,j=1, \ldots , n$ and $\nu , \mu \in I$ (where $\d_{ij}$ is the Kronecker delta): 
\begin{enumerate}
\item $[\d_{X_i},\d_{X_j}] = 0$,  $[\d_{Y_i},\d_{Y_j}] =0$, 
$ [ \d_{Y_i},\d_{X_j} ] =\d_{ij} \sum_{\nu \in I}\frac{\der}{\der t_\nu}(\der_i(a_i))\d_{t_\nu},  $\\
$ [ \d_{X_i},\d_{t_\nu} ] =-\der_i(t_\nu )\d_{X_i}-X_i \sum_{\mu \in I}\frac{\der}{\der t_\mu}(\der_i(t_\nu))\d_{t_\mu},  $\\
$ [ \d_{Y_i},\d_{t_\nu} ] =\der_i(t_\nu )\d_{Y_i}+Y_i \sum_{\mu \in I}\frac{\der}{\der t_\mu}(\der_i(t_\nu))\d_{t_\mu},  $
\item $[ \d_{X_i},t_\nu ] =-\der_i(t_\nu )X_i$, $[ \d_{Y_i},t_\nu ] =\der_i(t_\nu )Y_i$, 
$[ \d_{X_i},Y_j ] =-\d_{ij}\der_i(a_i)$, $[ \d_{Y_i},X_j ] =\d_{ij}\der_i(a_i)$,\\
$[ \d_{X_i},X_j ] =0$, $[ \d_{Y_i},Y_j ] =0$, $[ \d_{t_\nu},t_\mu ] =\{ t_\nu , t_\mu \}$,\\
$[ \d_{t_\nu},X_i ] =\der_i(t_\nu ) X_i$, $[ \d_{t_\nu},Y_i ] =-\der_i(t_\nu ) Y_i$, and
\item for all elements $s\in S$, $\sum_{\nu \in I}\frac{\der f_s}{\der t_\nu} \d_{t_\nu}
=0$ and  for all elements $i=1, \ldots , n$, $Y_i\d_{X_i}+X_i\d_{Y_i}-\sum_{\nu \in I} \frac{\der a_i}{\der t_\nu}\d_{t_\nu}=0$. 
\end{enumerate}
\end{theorem}

{\it Proof}. The theorem follows from
 Theorem \ref{23Jun19}.(2),  where the relations 1--3 are  the relations (a)--(c) of  Theorem \ref{23Jun19}.(2), respectively,  in the case of the GWPA $A$. $\Box$


\begin{proposition}\label{F23Jun19}
Let $A= K[H][X,Y; a, \der = b\der_H\}$ be a GWPA of degree $1$ where  $K[H]$ is a polynomial algebra in a single variable $H$ with trivial Poisson bracket, $a,b\in K[H]$ and $\der_H=\frac{d}{dH}=()'$. Then the algebra $\CU (A)$ is generated by the algebra $A=K[H,X,Y]/(XY-a)$ and the elements $\D_H$, $\d_X$ and $\d_Y$ subject to the defining relations given in statements 1-3: 
\begin{enumerate}
\item $[\d_Y,\d_X] = (ba')'\d_H$,  $[\d_X,\d_H] =-b\d_X-Xb'\d_H$, 
$ [ \d_Y,\d_H ] =b\d_Y+Yb'\d_H$, 
\item $[ \d_X,X] =0$, $[ \d_Y,Y] =0$, 
$[ \d_X,H] =-bX$, $[ \d_Y,H] =bY$,
$[ \d_X,Y ] =-ba'$, $[ \d_Y,X ] =ba'$, \\
$[ \d_H, H ] =0$, $[ \d_H,X ] =bX$, $[ \d_H,Y ] =bY$, and
\item $Y\d_X+X\d_Y=a'\d_H $. 
\end{enumerate}
\end{proposition}

{\it Proof}. The proposition follows from  Theorem \ref{E23Jun19}. The relations in statements 1--3 of the proposition are precisely the relations 1--3 of Theorem \ref{E23Jun19} for the GWPA $A$. $\Box$

\begin{corollary}\label{aF23Jun19}
Let the Poisson algebra $A =\bigotimes_{i=1}^n A_i$ be a tensor product of the GWPAs 
 $A_i= K[H_i][X_i,Y_i; a_i, \der_i = b_i\der_{H_i}\}$ from  Proposition \ref{F23Jun19}. Then $\CU (A) = \bigotimes_{i=1}^n\CU (A_i)$ and the defining relations for the algebra $\CU (A)$ are the union of the defining relations of all the algebra $\CU (A_i)$ from Proposition \ref{F23Jun19} and obvious commutation relations that come from the tensor product of algebras (i.e. $st=ts$).
\end{corollary}

{\it Proof}. \cite[Proposition 2.7]{GenDefRel-PUEA} states that the Poisson enveloping algebra of the tensor product of finitely many Poisson algebras is isomorphic to the tensor product of the corresponding Poisson enveloping algebras.  
 Now, the corollary follows from  Proposition \ref{E23Jun19}. $\Box$

\begin{corollary}\label{bF23Jun19}
Let $A= K[H][X,Y; a, \der = b\der_H\}$ be a GWPA of degree $1$ as in Proposition \ref{F23Jun19} such that $a\not \in K\cup \{ 0\}$, $A_{a'}$ be the localization of the algebra $A$ at the powers of the elements $a' = \frac{da}{dH}$.
   Then the algebra $\CU (A_{a'})$ is generated by the algebra $A_{a'}$ and the elements $\d_X$ and $\d_Y$ subject to the defining relations: 
\begin{enumerate}
\item $[\d_Y,\d_X] = \frac{(ba')'}{a'}(Y\d_X+X\d_Y )$,
\item $[ \d_X,X] =0$, $[ \d_Y,Y] =0$, 
$[ \d_X,H] =-bX$, $[ \d_Y,H] =bY$,
$[ \d_X,Y ] =-ba'$, $[ \d_Y,X ] =ba'$. 
\end{enumerate}
\end{corollary}

{\it Proof}. Proposition \ref{F23Jun19} gives generators and defining relations of the algebra $\CU (A)$. By  \cite[Theorem 2.10]{GenDefRel-PUEA},   $\CU (A_{a'})\simeq \CU (A)_{a'}$. Since $ \d_H=\frac{1}{a'} (Y\d_X+X\d_Y)$, the algebra $\CU (A_{a'})$ is generated by the algebra $A_{a'}$ and the elements $\d_X$ and $\d_Y$ subject to the defining relations in statements 1 and 2 of the corollary which are precisely the relations 1 and 2 of Proposition \ref{F23Jun19}, respectively (by direct computations, all the relations where $\d_H$ is involved in Proposition \ref{F23Jun19} are redundant).  $\Box$ \\


\section{A simplicity criterion for the Poisson  enveloping algebra of a  generalized Weyl Poisson algebra}\label{SIMPUEA}

{\bf Simplicity criterion for the algebra $P\CD (\CP )$ of Poisson differential operators on $\CP$.} An ideal $I$ of a Poisson algebra $\CP$ is called a {\bf Poisson ideal} if $\{\CP , I\}\subseteq I$. A Poisson algebra $\CP$ is called {\bf Poisson simple} if the ideals $0$ and $\CP$ are the only Poisson ideals of the Poisson algebra $\CP$. Let $\Der_K(\CP )$ be the Lie algebra of $K$-derivations of the (associative) algebra $\CP$. For each element $a\in \CP$, the derivation $\pad_a :=\{a, \cdot \}\in \Der_K(\CP)$ is called called the {\bf Hamiltonian vector field} associated with the element $a$.  Then the set of Hamiltonian vector fields 
$$\CH_\CP :=\{ \pad_a \, | \, a\in \CP\}$$ is a Lie subalgebra of the Lie algebra $\Der_K(\CP )$. \\

{\bf Definition.} Let $\CP$ be a Poisson algebra and $\CD (A)$ be the algebra of differential operators on $\CP$. The subalgebra of $\CD (\CP )$,

$$P\CD (\CP ):=\langle \CP , \CH_\CP \rangle,$$

is called the {\bf algebra of Poisson differential operators} of the Poisson algebra $\CP$.\\

In general, $P\CD (\CP )\neq \CD (\CP )$.
 Theorem \ref{X22Jul19} is a simplicity criterion for the algebra $P\CD (\CP )$ of Poisson differential operators on $\CP$.

\begin{theorem}\label{X22Jul19}
(\cite[Theorem 1.1]{GenDefRel-PUEA})  Let $\CP$ be a Poisson algebra over an arbitrary  field $K$.
Then the  following statements are equivalent:
\begin{enumerate}
\item The algebra $P\CD (\CP )$ is a simple algebra.
\item The Poisson algebra $\CP$ is a Poisson simple algebra. 
\end{enumerate}
\end{theorem}

{\bf Poisson simplicity criterion for  generalized Weyl Poisson   algebras.}
Let $D$ be a Poisson algebra and  $\der = (\der_1, \ldots , \der_n)\in \PDer_K(D)^n$. An ideal $I$ of $D$ is called $\der$-{\bf invariant} if $\der_i(I)\subseteq I$ for all $i=1, \ldots , n$. The set $$D^\der := \{ d\in D\, | \, \der_1(d)=0, \ldots , \der_n(d)=0\}$$ is called the {\bf algebra of $\der$-constants} of $D$, it is a subalgebra of $D$.  Theorem \ref{10Apr16} is  Poisson simplicity criterion for  generalized Weyl Poisson   algebras.

\begin{theorem}\label{10Apr16}
(\cite[Theorem 1.1]{Pois-GWA-Com}) Let $A=D[X,Y;a, \der\}$ be a GWPA of rank $n$. 
Then the  Poisson algebra $A$ is a simple Poisson algebra iff
\begin{enumerate}
\item the Poisson algebra  $D$ has no proper $\der$-invariant Poisson ideals,
\item for all $i=1, \ldots , n$, $Da_i+D\der_i(a_i) = D$, and
\item the algebra $\PZ (A)$ is a field $(\Leftrightarrow $   char$(K)=0$, $\PZ (D)^\der$ is a field and $D_{\alpha }=0$ for all $\alpha \in \Z^n\backslash \{ 0\} $ where  $D_\alpha := \{ \l \in D^\der\, | \, \pad_\l = \l \sum_{i=1}^n \alpha_i \der_i, \; \l \alpha_i\der_i(a_i)=0$ for $i=1, \ldots , n\}$,  \cite[ Proposition 1.2]{Pois-GWA-Com}).
\end{enumerate}
\end{theorem}

The next corollary follows from Theorem  \ref{X22Jul19} and Theorem \ref{10Apr16}.

\begin{corollary}\label{GWPA10Apr16}
 Let $A=D[X,Y;a, \der\}$ be a GWPA of rank $n$. 
Then the  following statements are equivalent:
\begin{enumerate}
\item The algebra $P\CD (A )$ is a simple algebra.
\item 
\begin{enumerate}
\item the Poisson algebra  $D$ has no proper $\der$-invariant Poisson ideals,
\item for all $i=1, \ldots , n$, $Da_i+D\der_i(a_i) = D$, and
\item the algebra $\PZ (A)$ is a field $(\Leftrightarrow $   char$(K)=0$, $\PZ (D)^\der$ is a field and $D_{\alpha }=0$ for all $\alpha \in \Z^n\backslash \{ 0\} $ where  $D_\alpha := \{ \l \in D^\der\, | \, \pad_\l = \l \sum_{i=1}^n \alpha_i \der_i, \; \l \alpha_i\der_i(a_i)=0$ for $i=1, \ldots , n\}$,  \cite[ Proposition 1.2]{Pois-GWA-Com}).
\end{enumerate}
\end{enumerate}
\end{corollary}

{\bf Simplicity criteria for the Poisson    enveloping algebra $\CU (\CP )$.} Let $\CP$ be a Poisson algebra. The algebra $\CP$ is a $\CD (\CP )$-module and hence $P\CD (\CP )$- and $\CU (\CP)$-module, and so 
 there is a natural algebra epimorphism (Theorem \ref{23Jun19}.(3)): 
 
\begin{equation}\label{UPDP}
\pi_\CP :\CU (\CP) \ra P\CD (\CP ),\;\; a\mapsto a, \;\; \d_b \mapsto \pad_b= \{ b, \cdot \}
\end{equation}
 where $a, b\in \CP$. 
 Theorem \ref{Y22Mar19} is a simplicity criteria for the Poisson    enveloping algebra $\CU (\CP )$.
 
\begin{theorem}\label{Y22Mar19}
(\cite[Theorem 1.2]{GenDefRel-PUEA})  Let $\CP$ be a Poisson algebra over an arbitrary  field $K$. Then the  following statements are equivalent:
\begin{enumerate}
\item The algebra $\CU (\CP )$ is a simple algebra.
\item The algebra $P\CD (\CP )$ is a simple algebra and $\CU (\CP )\simeq P\CD (\CP )$. 
\item  The Poisson algebra $\CP$ is a Poisson simple algebra and $\CP$ is a faithful left $\CU (\CP )$-module. 
\end{enumerate}
If one of the equivalent conditions holds then $\CU (\CP ) \simeq P\CD (\CP )$. 
\end{theorem}

Corollary \ref{GWPAY22Mar19}  is a simplicity criteria for the Poisson    enveloping algebra $\CU (A )$ where $A$ is a generalized Weyl Poisson algebra.   It  follows from Theorem \ref{10Apr16}  and Theorem \ref{Y22Mar19}.

\begin{corollary}\label{GWPAY22Mar19}
 Let $A=D[X,Y;a, \der\}$ be a GWPA of rank $n$. 
Then the  following statements are equivalent:
\begin{enumerate}
\item The algebra $\CU (A )$ is a simple algebra.
\item 
\begin{enumerate}
\item the Poisson algebra  $D$ has no proper $\der$-invariant Poisson ideals,
\item for all $i=1, \ldots , n$, $Da_i+D\der_i(a_i) = D$, 
\item the algebra $\PZ (A)$ is a field, and 

\item the algebra $A$ is  a faithful left $\CU (A)$-module. 
\end{enumerate}
\end{enumerate}
\end{corollary}

{\bf Generators and defining relations for the $\CA$-module $\Der_K(\CA )$ where $\CA$ is a regular domain of essentially finite type.} 
A localization of an affine commutative algebra is called an {\bf algebra of essentially finite type}. The following notation is fixed:  $P_n=K[x_1, \ldots , x_n]$ is a polynomial algebra over {\bf perfect} field 
$K$, $I=(f_1, \ldots , f_m)$ is a  prime but not a maximal ideal of $P_n$, $ \CA = S^{-1}(P_n/I)$ is a domain of essentially finite type and 
 $Q=Q(\CA )$ is its field of fractions,  $r= r\Big(\frac{\der  f_i}{\der x_j}\Big)$ is the rank (over $Q$) of the {\bf Jacobian matrix} $\Big(\frac{\der  f_i}{\der x_j}\Big)$ of $\CA$, $\ga_r$ is the {\bf Jacobian ideal} of the algebra $\CA$ which is
(by definition) generated by all the $r\times r$ minors of the
Jacobian matrix of  $\CA$ (the algebra $\CA$ is {\bf regular} iff $\ga_r=\CA$, it is the
{\bf Jacobian criterion of regularity}, \cite[Corollary 16.20]{Eisenbook}),
$\Omega_\CA$ is the {\bf module of  K\"{a}hler} differentials for the
algebra $\CA$. 

For $\i =(i_1, \ldots , i_r)$ such that $1\leq
i_1<\cdots <i_r\leq m$ and $\j =(j_1, \ldots , j_r)$ such that
$1\leq j_1<\cdots <j_r\leq n$, $\D
 (\i , \j )$ denotes the corresponding minor of the Jacobian matrix of $\CA$,  and the $\i$ (resp., $\j $) is called {\bf
non-singular} if $\D (\i , \j')\neq 0$ (resp., $\D (\i', \j )\neq
0$) for some $\j'$ (resp., $\i'$). We denote by $\mI_r$ (resp.,
$\mJ_r$) the set of all the non-singular $r$-tuples $\i$ (resp.,
$\j $).

Since $r$ is the rank of the Jacobian matrix of $\CA$, it is easy to
show that $\D (\i , \j )\neq 0$ iff $\i\in \mI_r$ and $\j\in
\mJ_r$, \cite[Lemma 2.1]{gendifreg}. We denote by $\mJ_{r+1}$ the set of all
$(r+1)$-tuples $\j =(j_1, \ldots , j_{r+1})$ such that $1\leq
j_1<\cdots <j_{r+1}\leq n$ and when deleting  some element, say
$j_\nu$, we have a non-singular $r$-tuple $(j_1, \ldots
,\widehat{j_\nu},\ldots , j_{r+1})\in \mJ_r$ where the hat over
 a symbol  means that the symbol is omitted from
the list. The set $\mJ_{r+1}$ is called the {\bf critical set} and
any element of it is called a {\bf critical singular}
$(r+1)$-tuple. $\Der_K(\CA)$ is the $\CA$-module of $K$-derivations of
the algebra $\CA$. 
 
 The next theorem gives a finite set of generators and a
finite set of  defining relations for the left $\CA$-module
$\Der_K(\CA)$ when $\CA$ is a  regular algebra.

\begin{theorem}\label{9bFeb05}
(\cite[Theorem 4.2]{gendifreg} if char$(K)=0$; \cite[Theorem 1.1]{gendifregcharp} if char$(K)=p>0$;) Let the algebra $\CA$ be a regular domain of essentially finite type over the perfect field $K$. Then the left $\CA$-module
$\Der_K(\CA)$ is generated by the derivations $\der_{\i , \j }$, $\i
\in \mI_r$, $\j \in \mJ_{r+1}$ where
\begin{eqnarray*}
  \der_{\i , \j }=\der_{i_1, \ldots , i_r; j_1, \ldots , j_{r+1}}:= \det
 \begin{pmatrix}
  \frac{\der f_{i_1}}{\der x_{j_1}} & \cdots &  \frac{\der f_{i_1}}{\der
  x_{j_{r+1}}}\\
  \vdots & \vdots & \vdots \\
  \frac{\der f_{i_r}}{\der x_{j_1}} & \cdots &  \frac{\der f_{i_r}}{\der
  x_{j_{r+1}}}\\
  \der_{j_1}& \cdots & \der_{j_{r+1}}\\
\end{pmatrix}
\end{eqnarray*}
that satisfy the following defining relations (as a left
$\CA$-module): 
\begin{equation}\label{Derel}
\D (\i , \j )\der_{\i', \j'}=\sum_{l=1}^s(-1)^{r+1+\nu_l}\D (\i';
j_1', \ldots , \widehat{j_{\nu_l}'}, \ldots , j_{r+1}')\der_{\i;\j
, j_{\nu_l}'}
\end{equation}
for all $\i, \i'\in \mI_r$, $\j=(j_1, \ldots , j_r)\in \mJ_r$, and
$\j'=(j_1', \ldots , j_{r+1}')\in \JJ_{r+1}$ where $\{ j_{\nu_1}',
\ldots , j_{\nu_s}'\}=\{ j_1', \ldots , j_{r+1}'\}\backslash \{
j_1, \ldots , j_r\}$.
\end{theorem}

Suppose that the algebra $\CA = S^{-1}(P_n/I)$ is a Poisson domain of essentially finite type. Let $d=d_\CA = r(\CC_\CA )$ be the rank of the $n\times n$ matrix $$\CC_\CA =(\{ x_i , x_j\} )\in M_n(\CA )$$ over
 the field $Q$. For each $l=1, \ldots , n$, let $$\ind_n(l) =\{ \i =(i_1, \ldots , i_l)\, | \, 1\leq i_1 <\cdots <i_l\leq n\}.$$ For  elements $\i = (i_1, \ldots , i_l)$ and $\j = (j_1, \ldots , j_l)$ of $\ind_n(l)$, let $\CC_\CA (\i , \j )=(\{ x_{i_\nu} , x_{j_\mu } \} )$ be the $l\times l$ submatrix of the matrix $\CC_\CA$. So, the rows (resp., the columns) of the matrix  $\CC_\CA (\i , \j )$ are indexed by $i_1, \ldots , i_l$ (resp., $j_1, \ldots , j_l$). The $(i_\nu , j_\mu )^{th}$ element of the matrix  $\CC_\CA (\i , \j )$ is $\{ x_{i_\nu} , x_{j_\mu }\}$. Let $M_{\CA , l} =\{ \CC_\CA (\i , \j ) \, | \, \i, \j \in \ind_n(l)\}$ be the set of all $l\times l$ submatrices of $\CC_\CA$ and 
 $$\CC_{\CA , l} =\{ \mu (\i , \j ) :=\det (\CC_\CA (\i , \j )) \, | \, \i, \j \in \ind_n(l)\}$$ be the set of all $l\times l$ minors of $\CC_\CA$. Let 
 \begin{eqnarray*}
 \mI (l)&=& \mI_\CA (l) = \{ \i \in \ind_n(l)\, | \, \mu (\i , \j ) \neq 0\;\; {\rm for \; some}\;\; \j \in \ind_n(l)\} \\
 \mJ (l)& = &\mJ_\CA (l) = \{ \j \in \ind_n(l)\, | \, \mu (\i , \j ) \neq 0\;\;  {\rm for \; some}\;\;\i \in \ind_n(l)\}
\end{eqnarray*}

{\bf The ideal $\kappa_\CA$ of the algebra $\CU (\CA )$ and its generators $\d_{\i , i_\nu; \j }$.} For each pair of elements $\i = (i_1, \ldots , i_d)$ and $\j = (j_1, \ldots , j_d)$ of $\mI_\CA (d)=\mJ_\CA (d)$ and each element $i_\nu\in \{ i_{d+1}, \ldots , i_n\}=\{ 1, \ldots , n\} \backslash \{ i_1, \ldots , i_d\}$, let us  consider the following elements of the algebra $\CU (\CA )$, \cite{GenDefRel-PUEA}, 
\begin{eqnarray}\label{dinnj}
 \d_{\i , i_\nu; \j }:= \det
 \begin{pmatrix}
 \{ x_{i_1}, x_{j_1} \} & \ldots &   \{ x_{i_1}, x_{j_d}\} & \d_{i_1}\\
  \vdots & \vdots & \vdots  & \vdots\\
   \{ x_{i_d}, x_{j_1} \} & \ldots &   \{ x_{i_d}, x_{j_d} \}& \d_{i_d}\\
  \{ x_{i_\nu}, x_{j_1} \} & \ldots &   \{ x_{i_\nu}, x_{j_d} \}& \d_{i_\nu}\\
\end{pmatrix}= \mu (\i , \j ) \d_{i_\nu} + \sum_{s=1}^d(-1)^{s+d+1}\mu (i_1, \ldots ,\widehat{i_s}, \ldots , i_d, i_\nu;\j )\d_{i_s} 
\end{eqnarray}
where $\d_i:=dx_i\in \O_\CA\subseteq \CU (\CA )$.\\

{\it Definition, \cite{GenDefRel-PUEA}.} Let $\kappa_\CA$ be an ideal of the algebra $\CU (\CA )$ which is generated by the finite set of   elements  $\d_{\i , i_\nu; \j }\in \O_\CA$ where $\i , \j \in \mI_\CA (d)$, see (\ref{dinnj}). Then  $\kappa_\CA\subseteq \ker (\pi_\CA )$, \cite{GenDefRel-PUEA}. 
\\

 {\bf Criteria for $\ker (\pi_\CA )=0$.}
 In the case when the Poisson algebra $\CA$ is a {\em regular} domain of essentially finite type, Theorem \ref{3Aug19} is an efficient explicit criterion for $\ker (\pi_\CA )=0$, i.e. for the epimorphism $\pi_\CA : \CU (\CA ) \ra P\CD (\CA )$ to be an isomorphism.

 Since $\Der_K(\CA ) \simeq \Hom_\CA (\O_\CA , \CA )$, there is a {\em pairing} of left $\CA$-modules  (which is  an $\CA$-bilinear map):
$$
\Der_K(\CA ) \times \O_\CA \ra \CA , \;\; (\der , \o ) \mapsto (\der, \o ) := \der (\o ).
$$

\begin{theorem}\label{3Aug19}
(\cite[Theorem 1.9]{GenDefRel-PUEA}) 
Let a Poisson algebra $\CA= S^{-1}(P_n/I)$ be a regular domain of essentially finite type over the field $K$ of characteristic zero, $d=r(\CC_\CA )$ and $r=r\Big( \frac{\der f_i}{\der x_j}\Big)$. Then the following statements are equivalent (the derivations $\der_{\i ; \j , j_\nu}$ of $\CA$ are defined in Theorem \ref{9bFeb05}):
\begin{enumerate}
\item $\ker (\pi_\CA )=0$ $(\Leftrightarrow \pi_\CA : \CU (\CA ) \simeq P\CD (\CA ))$.  
\item $\kappa_\CA =0$. 
\item $d=n-r$ and $(\der_{\i ; \j , j_\nu}, \d_{\i' ; \j' , j_\mu'})=0$ for all elements $\i \in \mI_r$, $\j\in \mJ_r$, $\i'\in \mI_\CA (d)$, $\j'\in \mJ_\CA (d)$, $\nu = r+1, \ldots , n$ and $\mu = d+1, \ldots , n$ where for $\j = (j_1, \ldots , j_r)$ and $\i' = (i_1', \ldots , i'_d)$, $\{ j_{r+1}, \ldots , j_n\} =\{ 1,\ldots , n\} \backslash \{ j_1, \ldots , j_r\}$ and $\{ i'_{d+1}, \ldots , i'_n\} =\{ 1, \ldots , n\} \backslash \{ i_1', \ldots , i_d'\}$. 
\end{enumerate}
\end{theorem}


   $\GK$ stands for the {\bf Gelfand-Kirillov dimension} and  $\Sym_\CA (\O_\CA )$  is the {\bf symmetric algebra} of the $\CA$-module  $\O_\CA$ of  K\"{a}hler differentials of the  algebra $\CA$.

In the case when the Poisson algebra $\CP = \CA$ is an algebra  of essentially finite type over a field of characteristic zero,  Theorem \ref{Y22Mar19} can be strengthen, see Theorem \ref{7Aug19}.

\begin{theorem}\label{7Aug19}
(\cite[Theorem 1.3]{GenDefRel-PUEA}) Let  a Poisson algebra $\CA$ be an algebra of essentially finite type over the field $K$ of characteristic zero.
 Then the  following statements are equivalent:
\begin{enumerate}
\item The algebra  $\CU (\CA )$ is a simple algebra.
\item The algebra $P\CD (\CA )$ is a simple algebra and  one of the equivalent conditions of Theorem \ref{3Aug19} holds.
\item  The algebra $\CA$ is  Poisson simple  and one of the equivalent conditions of Theorem \ref{3Aug19} holds. 
\end{enumerate}
If one of the equivalent conditions holds then the algebra $\CA = S^{-1}(P_n/I)$ is a regular, Poisson simple  domain of essentially finite type over the  field $K$ of characteristic zero, the algebra epimorphism 
$\pi_\CA :  \CU (\CA) \ra  P\CD (\CA )$ is an isomorphism, 
 $d=n-r$ where  $d= r(\CC_\CA )$ and $r= r\Big(\frac{\der  f_i}{\der x_j}\Big)$, and the algebra $\CU (\CA )$ is a simple Noetherian domain with  $$\GK \, \CU (\CA ) = \GK \, P\CD (\CA) = \GK \, {\rm gr}\, \CU (\CA ) = \GK \, \Sym_\CA (\O_\CA )= 2\GK (\CA ) = 2(n-r).$$ 
\end{theorem}

{\bf A simplicity criterion for the PEA $\CU (A )$ of a GWPA $A$.} Let $A =D[X,Y; \der , a\}$ be a GWPA of degree $N$; a Poisson algebra $D = \CA =S^{-1}(P_n/I)$ be a regular domain of essentially finite type; $I=(f_1, \ldots , f_m)$; $r= r\Big(\frac{\der  f_i}{\der x_j}\Big)$ be  the rank of the Jacobian matrix $\Big(\frac{\der  f_i}{\der x_j}\Big)$ of the algebra $D$; $d= r(\CC_A )$ be the rank of the $(n+2N)\times (n+2N)$ matrix

$$ \CC_A = (\{ s_i, s_j\})\;\; {\rm where}\;\; s_i, s_j\in \{ x_1, \ldots , x_n, X_1, \ldots , X_N, Y_1, \ldots , Y_N\}.$$
 
Theorem \ref{GWPA7Aug19} is a simplicity criterion for the Poisson enveloping algebra $\CU (A )$ of a GWPA $A$.

\begin{theorem}\label{GWPA7Aug19}
Let  a Poisson algebra $D = S^{-1}(P_n/I)$ be a regular domain of essentially finite type  over the field $K$ of characteristic zero  and $r= r\Big(\frac{\der  f_i}{\der x_j}\Big)$ be the rank of the Jacobian matrix of $D$, $A =D[X,Y; \der , a\}$ be a GWPA of degree $N$ and  
$d= r(\CC_A )$.  Then the  following statements are equivalent:
\begin{enumerate}
\item The algebra $\CU (A )$ is a simple algebra.
\item The algebra $P\CD (A )$ is a simple algebra
and   one of the equivalent conditions of Theorem \ref{3Aug19} holds.
\item  The Poisson algebra $A$ is a Poisson simple algebra  
and   one of the equivalent conditions of Theorem \ref{3Aug19} holds. 
\item 
\begin{enumerate}
\item the Poisson algebra  $D$ has no proper $\der$-invariant Poisson ideals,
\item for all $i=1, \ldots , n$, $Da_i+D\der_i(a_i) = D$, and
\item the algebra $\PZ (A)$ is a field $(\Leftrightarrow $   char$(K)=0$, $\PZ (D)^\der$ is a field and $D_{\alpha }=0$ for all $\alpha \in \Z^n\backslash \{ 0\} $ where  $D_\alpha := \{ \l \in D^\der\, | \, \pad_\l = \l \sum_{i=1}^n \alpha_i \der_i, \; \l \alpha_i\der_i(a_i)=0$ for $i=1, \ldots , n\}$,  \cite[ Proposition 1.2]{Pois-GWA-Com}).
\item $d=n-r+N$ and $(\der_{\i ; \j , j_\nu}, \d_{\i' ; \j' , j_\mu'})=0$ for all elements $\i \in \mI_r$, $\j\in \mJ_r$, $\i'\in \mI_\CA (d)$, $\j'\in \mJ_\CA (d)$, $\nu = r+1, \ldots , n$ and $\mu = d+1, \ldots , n$ where for $\j = (j_1, \ldots , j_r)$ and $\i' = (i_1', \ldots , i'_d)$, $\{ j_{r+1}, \ldots , j_n\} =\{ 1,\ldots , n\} \backslash \{ j_1, \ldots , j_r\}$ and $\{ i'_{d+1}, \ldots , i'_n\} =\{ 1, \ldots , n\} \backslash \{ i_1', \ldots , i_d'\}$. 
\end{enumerate}
\end{enumerate}
If one of the equivalent conditions holds then 
$\CU (A) \simeq P\CD (A )$, the algebra $A$ is a regular, Poisson simple domain (in particular, all $a_i\neq 0$), $d=n-r+N$ and the algebra $\CU (A )$ is a simple Noetherian domain with  $$\GK \, \CU (A ) = \GK \, P\CD (A) =  2\GK (A ) = 2(n-r+N).$$ 
\end{theorem}

{\it Proof}. The GWPA $A$ is a domain iff  the algebra $D$  and all the  elements $a_i$ are not equal to zero.
Notice that $\GK (A)=\GK (D)+N=n-r+N$. Now, the theorem follows from Theorem \ref{7Aug19}, Theorem  \ref{10Apr16} and Theorem \ref{3Aug19}. $\Box$\\

So, Theorem \ref{GWPA7Aug19}.(4) is an efficient tool in proving or disproving simplicity of the algebra $\CU (A )$.\\

 {\bf The Gelfand-Kirillov dimension of the algebras $\CU (A )$,    ${\rm gr}\, \CU (A )$ and $ \Sym_A (\O_A )$.} 
 
\begin{theorem}\label{AAA29Jul19}
(\cite[Theorem 1.4]{GenDefRel-PUEA}) Let a Poisson algebra $\CA = S^{-1} (P_n/I)$ be a  domain of essentially finite type over a perfect field $K$ where  $I=(f_1, \ldots , f_m)$ is a prime   ideal of $P_n$,  and $r=r\Big(\frac{\der f_i}{\der x_j}\Big)$ be the rank of the Jacobian matrix $\Big(\frac{\der f_i}{\der x_j}\Big)$ .  Then the algebra $\CU (\CA )$ is a Noetherian algebra with $$\GK \, \CU (\CA ) = \GK \, {\rm gr}\, \CU (\CA ) = \GK \, \Sym_\CA (\O_\CA ) = 2\GK (\CA ) = 2(n-r).$$
\end{theorem}

By Theorem  \ref{E23Jun19}, the Poisson enveloping algebra $\CU (A)$ of  the GWPA 
$A= D[X,Y; a, \der \}$, where  $D = S^{-1}(P_n/I)$, admits a filtration by the total degree of the elements $\d_{x_1}, \ldots, \d_{x_n}$, $\d_{X_1}, \ldots, \d_{X_n}$, $ \d_{Y_1}, \ldots, \d_{Y_n} $ (nonzero  elements of the algebra $A$ have degree zero). The associated graded algebra of $\CU (A)$ is denoted by 
 ${\rm gr}\, \CU (A )$. In the case when $\CU (A)=P\CD (A)$ and char$(K)=0$, the filtration coincides with the order filtration on the algebra of Poisson differential operators $P\CD (A)$ on $A$.

\begin{corollary}\label{GAAA29Jul19}
Let  a Poisson algebra $D = S^{-1}(P_n/I)$ be a  domain of essentially finite type over a perfect field $K$, $I=(f_1, \ldots , f_m)$, $r= r\Big(\frac{\der  f_i}{\der x_j}\Big)$ be the rank of the Jacobian matrix $\Big(\frac{\der f_i}{\der x_j}\Big)$ of $D$, $A =D[X,Y; \der , a\}$ be a GWPA of degree $N$ with all $a_i\neq 0$ and  $d= r(\CC_A )$. Then the algebra $\CU (A )$ is a Noetherian algebra with $$\GK \, \CU (A ) = \GK \, {\rm gr}\, \CU (A ) = \GK \, \Sym_A (\O_A ) = 2\GK (A ) = 2(n-r+N).$$
\end{corollary}

{\it Proof.}  Since $\GK (A)=\GK (D)+N=n-r+N$, the corollary follows from Theorem \ref{AAA29Jul19}. $\Box$\\

{\bf The Gelfand-Kirillov dimension  of 
the algebra $P\CD (A )$ of Poisson differential operators on $A $.} Proposition \ref{BB26Jul19} gives the exact figure for the Gelfand-Kirillov dimension  of 
the algebra $P\CD (\CA )$ where the Poisson algebra $\CA$ is a domain of essentially finite type over the field $K$ of characteristic zero.

\begin{proposition}\label{BB26Jul19}
(\cite[Theorem 1.5]{GenDefRel-PUEA})  Let a Poisson algebra $\CA$ be a domain of essentially finite type over the field $K$ of characteristic zero and  $r$ be the rank of Jacobian matrix of $\CA$  and $d=r(\CC_\CA )$. Then 
 $$\GK (P\CD (\CA )) = \GK (\CA ) +d= n-r+d .$$
\end{proposition}

\begin{corollary}\label{GBB26Jul19}
 
Let  a Poisson algebra $D = S^{-1}(P_n/I)$ be a  domain of essentially finite type over the  field $K$ of characteristic zero, $I=(f_1, \ldots , f_m)$, $r= r\Big(\frac{\der  f_i}{\der x_j}\Big)$ be the rank of the Jacobian matrix $\Big(\frac{\der f_i}{\der x_j}\Big)$ of $D$, $A =D[X,Y; \der , a\}$ be a GWPA of degree $N$ with all $a_i\neq 0$ and  $d= r(\CC_A )$. Then 
 $$\GK (P\CD (A )) = \GK (A ) +d= n-r+N+d .$$
\end{corollary}

{\it Proof.}  The GWPA $A$ is a domain since the algebra $D$ is so and all the  elements $a_i$ are not equal to zero.  Since $\GK (A)=\GK (D)+N=n-r+N$, the corollary follows from Theorem \ref{BB26Jul19}. $\Box$\\

 {\bf The algebra $\CU (A )$ is a domain when $A$ is a regular domain of essentially finite type.}  Theorem \ref{29Jul19} states that the algebra $\CU (\CA )$ is a domain provided the algebra $\CA$ is a  regular domain  of essentially finite type. 

\begin{theorem}\label{29Jul19}
(\cite[Theorem 4.4]{GenDefRel-PUEA}) Let a Poisson algebra $\CA = S^{-1} (P_n/I)$ be a regular domain of essentially finite type over the perfect field $K$  where  $I=(f_1, \ldots , f_m)$ is a prime but not maximal  ideal of $P_n$ and $r=r\Big(\frac{\der f_i}{\der x_j}\Big)$ is the rank of the Jacobian matrix $\Big(\frac{\der f_i}{\der x_j}\Big)$ over the field of fractions of the domain $P_n/I$.  Then the algebra $\CU (\CA )$ is a Noetherian domain with $\GK \, \CU (\CA ) = \GK \, {\rm gr}\, \CU (\CA ) = \GK \, \Sym_\CA (\O_\CA ) = 2\GK (\CA ) = 2(n-r)$. 
\end{theorem}

\begin{corollary}\label{G29Jul19}

 Let  a Poisson algebra $D = S^{-1}(P_n/I)$ be a  domain of essentially finite type over the perfect field $K$, $I=(f_1, \ldots , f_m)$, $r= r\Big(\frac{\der  f_i}{\der x_j}\Big)$ be the rank of the Jacobian matrix $\Big(\frac{\der f_i}{\der x_j}\Big)$ of $D$, $A =D[X,Y; \der , a\}$ be a regular GWPA of degree $N$ with all $a_i\neq 0$ and  $d= r(\CC_A )$. Then the algebra $\CU (A )$ is a Noetherian domain with $\GK \, \CU (A ) = \GK \, {\rm gr}\, \CU (A ) = \GK \, \Sym_A (\O_A ) = 2\GK (A ) = 2(n-r+N)$. 
\end{corollary}

{\it Proof.}  Since $\GK (A)=\GK (D)+N=n-r+N$ and the algebra $A$ is a domain (since $D$ is so and all the elements $a_i\neq 0$), the corollary follows from Theorem \ref{29Jul19}. $\Box$\\

{\bf Data Availability.} Data sharing is not applicable to this article as no new data were created or analyzed in this study.

{\small

Department of Pure Mathematics

University of Sheffield

Hicks Building

Sheffield S3 7RH

UK

email: v.bavula@sheffield.ac.uk
}

\end{document}